\documentclass[12pt]{article}
\usepackage{latexsym}
\font\smallit=cmti10

\begin{document} 
\begin{center}
{\bf \large Parity-Alternate Permutations and Signed Eulerian Numbers}
\vskip 20pt
{\bf Shinji Tanimoto}\\
{\smallit Department of Mathematics, Kochi Joshi University, Kochi 780-8515, Japan}\\
{\tt tanimoto@cc.kochi-wu.ac.jp}\\ 
\end{center}
\begin{abstract}
In order to study signed Eulerian numbers, we introduce permutations of a particular type, 
called parity-alternate permutations, 
because they take even and odd entries alternately. 
The objective of this paper is twofold.
The first is to derive several properties of those permutations, by subdividing them 
into even and odd ones. 
The second is to discuss relationships between those and signed Eulerian numbers. 
Divisibility properties by prime powers are also deduced 
for signed Eulerian numbers and several related numbers. \\
\end{abstract}

\begin{center}
{\bf 1. Introduction}
\end{center}
Let $n$ be a positive integer and $A=a_1a_2\cdots a_n$ a permutation of $[n]=\{1,2, \ldots, n\}$.
An {\it ascent} of $A$ is an adjacent pair such that $a_i < a_{i+1}$ for some $i$ ($1\le i \le n-1$).  
Let $E(n,k)$ be the set of all permutations of $[n]$ 
with exactly $k$ ascents, where $0 \le k \le n-1$. Its cardinality is the classical
{\it Eulerian number};
\[
    A_{n,k}= |E(n,k)|, 
\]
whose properties and identities can be found in [2-6]. In particular, the well-known recurrence 
relation holds:
\begin{eqnarray}
A_{n,k} =  (n-k)A_{n-1,k-1} + (k+1)A_{n-1,k}.
\end{eqnarray} 
\indent
An {\it inversion} of a permutation $A=a_1a_2\cdots a_n$ is a 
pair $(i, j)$ such that $1 \le i < j \le n$ and $a_i > a_j$. Let us denote by ${\rm inv}(A)$   
the number of inversions in $A$. By $E_{\rm e}(n,k)$ and $E_{\rm o}(n,k)$ we denote
the subsets of all permutations in $E(n,k)$ 
that have, respectively, even and odd numbers of inversions, and 
their cardinalities by
\[
     B_{n,k} = |E_{\rm e}(n,k)|~~  {\rm and} ~~  C_{n,k} = |E_{\rm o}(n,k)|.
\]   
Obviously we have $A_{n,k}=B_{n,k}+C_{n,k}$, while the differences
\[
        D_{n,k}=B_{n,k}- C_{n,k}
	\]
were called {\it signed Eulerian numbers} in [1, 7], in which the descent number 
was considered instead of the ascent number. Therefore, the identities for $D_{n,k}$ 
presented here correspond to those in [1, 7] that are obtained by replacing $k$ with $n-k-1$. \\
\indent
In [1, 10] it was proved that the recurrence relation for $D_{n,k}$ is 
\begin{eqnarray}
     D_{n,k} = 
        \left\{\begin{array}{ll}
                    (n-k)D_{n-1,k-1}+ (k+1)D_{n-1,k}, & \mbox{if  $n$ is odd}, \\
                    D_{n-1,k-1} - D_{n-1,k}, &  \mbox{if  $n$ is even}. \\
                    \end{array}  \right.
\end{eqnarray}
The values of $D_{n,k}$ for small $n$ are given below for reference. Those for $B_{n,k}$ and $C_{n,k}$
can be found in [10]. \\
\\
\begin{tabular}{l|rrrrrrrrrrr}
         \noalign{\hrule height 0.8pt}
        $D_{n,k}$ & 0 & 1 & 2 & 3 & 4 & 5 & 6 & 7 & ~8 & ~~9  \\
       \hline
         $n = 2$ & $-1$ & 1 \\
         $n = 3$ & $-1$ & 0 & 1  \\
         $n = 4$ & 1 & $-1$ & $-1$ & 1 \\
         $n = 5$ & 1 & 2 & $-6$ & 2 & 1 \\
         $n = 6$ & $-1$ & $-1$ & 8 & $-8$ & 1 & 1 \\
         $n = 7$ & $-1$ & $-8$ & 19 & 0 & $-19$ & 8 & 1 \\ 
         $n = 8$ & 1 & 7 & $-27$ & 19 & 19 & $-27$ & 7 & 1 \\ 
         $n = 9$ & 1 & 22 & $-32$ & $-86$ & 190 & $-86$ & $-32$ & 22 &1 \\ 
	 $n = 10$ & $-1$ & $-21$ & 54 & 54 & $-276$ & 276 & $-54$ & $-54$ & 21 & 1\\ 
\noalign{\hrule height 0.8pt}
\end{tabular}
\\
\\
\\
\indent
The objective of this paper is to investigate permutations of a particular type, 
which we call parity-alternate, because
those take even and odd integers alternately, such as 236145 or 5274163. Utilizing such permutations,
we study signed Eulerian numbers and several related numbers. \\
\indent
In Section 2 parity-alternate permutations are introduced and their fundamental properties are
studied, by subdividing them into even and odd ones.
In Sections 3 and 4 we discuss further properties concerning those permutations using the operator
introduced in [8]. In particular, close relationships to signed Eulerian numbers are exhibited.
In Section 5 divisibility properties of signed Eulerian numbers and their related ones 
are considered as in [9] for classical Eulerian numbers. Moreover,
a similar relation to (1) will be deduced for the cardinalities for parity-alternate permutations.\\
\begin{center}
{\bf 2. Parity-Alternate Permutations}
\end{center}
A permutation will be called {\it parity-alternate permutation} ({\it PAP})
if its entries take even and odd integers alternately as
$436125$, $563412$ or $7216345$, for example. Note that, when $n$ is an odd integer, 
odd entries must appear at both ends of PAPs. Hence the total number of PAPs of $[n]$ is equal to 
\[
                    2\left(\left(\frac{n}{2}\right)!\right)^2,
\]
when $n$ is even,  and
\[
            \left(\frac{n+1}{2}\right)!\left(\frac{n-1}{2}\right)! =
       \frac{n+1}{2} \left(\left(\frac{n-1}{2}\right)!\right)^2, 
\]
when $n$ is odd.\\
\indent
We denote by $\Xi(n,k)$ the set of all PAPs in $[n]$ with exactly $k$ ascents and
divide them into two according to the parity of permutations:
\[
\Xi_{\rm e}(n,k) = \Xi(n,k) \cap E_{\rm e}(n,k),~~
\Xi_{\rm o}(n,k) = \Xi(n,k) \cap E_{\rm o}(n,k).
\]
\indent
The objective of this paper is to investigate their cardinalities 
\[
S_{n,k} = |\Xi(n,k)|, ~~ P_{n,k} = |\Xi_{\rm e}(n,k)|, ~~  Q_{n,k} = |\Xi_{\rm o}(n,k)|,
\]
together with the difference 
\[
 R_{n,k} = P_{n,k} - Q_{n,k}.
 \]
The last one is closely related to signed Eulerian numbers. In particular, we will show that the equality
\begin{eqnarray*}
  D_{n,k}=R_{n,k}
\end{eqnarray*}  
holds.
Obviously we have $S_{n,k} = P_{n,k} + Q_{n,k}$. \\
\indent
It is easy to see that the number of even PAPs is equal to that of odd ones.
When one interchanges 1 with 3 in PAPs, the resulting permutations are also
parity-alternate with opposite parity and this operation is a bijection on the set of such permutations. 
Hence the total number of even PAPs of $[n]$ is equal to 
\[
                    \left(\left(\frac{n}{2}\right)!\right)^2,
\]
when $n$ is even,  and
\[
            \frac{1}{2}\left(\frac{n+1}{2}\right)!\left(\frac{n-1}{2}\right)! =
       \frac{n+1}{4} \left(\left(\frac{n-1}{2}\right)!\right)^2, 
\]
when $n$ is odd. \\
\indent
For a permutation $A = a_1a_2 \cdots a_n$ we define its reflection 
by $A^{\ast}= a_n  \cdots a_2a_1$. The fact that
$A \in \Xi(n,k)$ if and only if $A^{\ast} \in \Xi(n,n-k-1)$
yields the symmetry property:
\begin{eqnarray}
S_{n,k}=S_{n,n-k-1}. 
\end{eqnarray}
Moreover, using reflected permutations and the parity of $n(n-1)/2$, the following symmetry 
properties between $P_{n,k}$ and $Q_{n,k}$ are easily checked. \\
\indent
(i) $n \equiv$ 0  or 1 ($\bmod$ 4).  In this case, $A \in E_{\rm e}(n,k)$ if and only if
$A^{\ast} \in E_{\rm e}(n,n-k-1)$, and 
$A \in E_{\rm o}(n,k)$ if and only if
$A^{\ast} \in E_{\rm o}(n,n-k-1)$, so we have
\[
     P_{n,k} = P_{n,n-k-1}~~ {\rm and} ~~  Q_{n,k} = Q_{n,n-k-1}.
\]
 \indent
 (ii) $n \equiv$ 2 or 3 ($\bmod$ 4). In this case, 
$A \in E_{\rm e}(n,k)$ if and only if
$A^{\ast} \in E_{\rm o}(n,n-k-1)$, and 
$A \in E_{\rm o}(n,k)$ if and only if
$A^{\ast} \in E_{\rm e}(n,n-k-1)$, so we have
\[
     P_{n,k} = Q_{n,n-k-1}~~ {\rm and} ~~ Q_{n,k} = P_{n,n-k-1}.
\]
\indent
The values of $S_{n,k}$, $P_{n,k}$ and $Q_{n,k}$ for small $n$ are presented below. \\
\\
\begin{tabular}{l|rrrrrrrrrrr}
         \noalign{\hrule height 0.8pt}
        $S_{n,k}$ & 0 & 1 & 2 & 3 & 4 & 5 & 6 & 7 & 8 & ~9 \\
       \hline
         $n = 2$ & 1 & 1 \\
         $n = 3$ & 1 & 0 & 1  \\
         $n = 4$ & 1 & 3 & 3 & 1 \\
         $n = 5$ & 1 & 2 & 6 & 2 & 1 \\
         $n = 6$ & 1 & 9 & 26 & 26 & 9 & 1 \\
         $n = 7$ & 1 & 8 & 39 & 48 & 39 & 8 & 1 \\ 
         $n = 8$ & 1 & 23 & 165 & 387 & 387 & 165 & 23 & 1 \\ 
         $n= 9$ & 1 & 22 & 228 & 674 & 1030 & 674 & 228 & 22 &~1 \\
         $n= 10$ & 1 & 53 & 860 & 4292 & 9194 & 9194 & 4292 & 860 & 53 &~1 \\
\noalign{\hrule height 0.8pt}
\end{tabular}
\\
\\
\\
\begin{tabular}{l|rrrrrrrrrrr}
         \noalign{\hrule height 0.8pt}
        $P_{n,k}$ & 0 & 1 & 2 & 3 & 4 & 5 & 6 & ~7 & ~8 & ~9 \\
       \hline
         $n = 2$ & 0 & 1 \\
         $n = 3$ & 0 & 0 & 1  \\
         $n = 4$ & 1 & 1 & 1 & 1 \\
         $n = 5$ & 1 & 2 & 0 & 2 & 1 \\
         $n = 6$ & 0 & 4 & 17 & 9 & 5 & 1 \\
         $n = 7$ & 0 & 0 & 29 & 24 & 10 & 8 & 1 \\ 
         $n = 8$ & 1 & 15 & 69 & 203 & 203 & 69 & 15 & 1 \\ 
         $n= 9$ & 1 & 22 & 98 & 294 & 610 & 294 & 98 & 22 &~1 \\
	 $n= 10$ & 0 & 16 & 457 & 2173 & 4459 & 4735 & 2119 & 403 & 37 & 1\\
\noalign{\hrule height 0.8pt}
\end{tabular}
\\
\\
\\
\begin{tabular}{l|rrrrrrrrrrr}
         \noalign{\hrule height 0.8pt}
        $Q_{n,k}$ & 0 & ~1 & 2 & 3 & 4 & 5 & 6 & ~7 & ~8 & ~9  \\
       \hline
         $n = 2$ & 1 & 0 \\
         $n = 3$ & 1 & 0 & 0  \\
         $n = 4$ & 0 & 2 & 2 & 0 \\
         $n = 5$ & 0 & 0 & 6 & 0 & 0 \\
         $n = 6$ & 1 & 5 & 9 & 17 & 4 & 0 \\
         $n = 7$ & 1 & 8 & 10 & 24 & 29 & 0 & 0 \\ 
         $n = 8$ & 0 & 8 & 96 & 184 & 184 & 96 & 8 & 0 \\ 
         $n= 9$ & 0 & 0 & 130 & 380 & 420 & 380 & 130 & 0 & 0 \\
	 $n= 10$ & 1 & 37 & 403 & 2119 & 4735 & 4459 & 2173 & 457 & 16 & 0\\ 
\noalign{\hrule height 0.8pt}
\end{tabular}
\\
\\
\begin{center}
{\bf 3. Periods and orbits in {\boldmath $\Xi(n,k)$}}
\end{center}
In order to study properties of $S_{n,k}$, $P_{n,k}$ and $Q_{n,k}$,
we review the operator on permutations of $[n]$ that was introduced in [8]. 
The operator $\sigma$ is defined by adding one to 
all entries of a permutation and by changing $n+1$ into one. 
However, when $n$ appears at either end of a permutation, it is removed
and one is put at the other end. That is, for a permutation $a_1a_2\cdots a_n$ with $a_i = n$ 
for some $i$ $(2 \le i\le n-1)$, we have \\
\indent
(i) $\sigma(a_1a_2\cdots a_n) = b_1b_2\cdots b_n$, \\
where $b_i = a_i +1$ for all $i$ $(1 \le i\le n)$ and $n+1$ is replaced by one. 
And, for a permutation $a_1a_2\cdots a_{n-1}$ of $[n-1]$, we have:
\vspace{-0.5cm}
\begin{itemize}
\item[(ii)] $\sigma(a_1a_2\cdots a_{n-1}n) = 1b_1b_2\cdots b_{n-1}$; 
\item[(iii)] $\sigma(na_1a_2\cdots a_{n-1}) = b_1b_2\cdots b_{n-1}1$, 
\end{itemize}
\vspace{-0.5cm}
where $b_i = a_i +1$ for all $i$ $(1 \le i\le n-1)$. We denote by $\sigma^{\ell}A$
the repeated $\ell$ applications of $\sigma$ to a permutation $A$. \\
\indent
Let us denote by $E^{\lhd}(n,k)$ the set of all permutations
$a_1a_2\cdots a_n$ in $E(n,k)$ with $a_1 < a_n$. 
It is obvious that the operator preserves the number of ascents in a permutation and that the relation holds:
\begin{eqnarray*}
                    \sigma A \in E^{\lhd}(n,k)  ~~\mbox{if and only if}~~  A \in E^{\lhd}(n,k).                
\end{eqnarray*}
\indent
Let us observe the number of inversions of a permutation when $\sigma$ is applied. 
When $n$ appears at either end of a permutation $A=a_{1}a_{2}\cdots a_{n}$
as in (ii) or (iii), it is evident that 
\[
{\rm inv}(\sigma A) = {\rm inv}(A).
\]
Next consider the case (i). When $a_{i}=n$ for some $i$ ($2\le i \le n-1$), we get
$\sigma (a_1a_{2}\cdots a_{n}) =b_1b_2 \cdots b_n$, 
where $b_{i}= 1$ is at the $i$th position. In this case, $n-i$ inversions 
$(i, i+1), \ldots, (i, n)$ of $A$ vanish and, in turn, 
$i-1$ inversions $(1, i), \ldots, (i-1, i)$ of $\sigma A$ occur. Hence
the difference between the numbers of inversions is 
\begin{eqnarray}
 {\rm inv}(\sigma A) - {\rm inv}(A) = (i-1) - (n-i) = 2i - (n + 1).
\end{eqnarray}
\indent
This means that, when $n$ is even, each application of the operator changes the parity 
of permutations as long as $n$ remains in the interior of permutations. 
If $n$ is odd, however, the operator $\sigma$ also preserves the parity of all 
permutations of $[n]$.\\
\indent
Let us put
\[ 
\Xi^{\lhd}(n,k) = \Xi(n,k) \cap E^{\lhd}(n,k),
\]
and 
\[ 
\Xi^{\lhd}_{\rm e}(n,k) = \Xi^{\lhd}(n,k) \cap E_{\rm e}(n,k), ~~ 
\Xi^{\lhd}_{\rm o}(n,k) = \Xi^{\lhd}(n,k) \cap E_{\rm o}(n,k).
\]
When $n$ is odd, for a PAP 
$A=a_1a_2\cdots a_n$ such that $a_i = n$ for some $i$ $(2 \le i\le n-1)$, 
$\sigma A$ is not parity-alternate, since $(i-1)$th, $i$th and $(i+1)$th entries of $\sigma A$ are all odd,
as is seen from (i). \\
\indent
Note that in the case of even $n$, however, we have the relations:
\begin{eqnarray*}
               &    \sigma A \in \Xi(n,k)   ~~\mbox{if and only if }~~ A \in \Xi(n,k); \\
	     &    \sigma A \in \Xi^{\lhd}(n,k)  ~~\mbox{if and only if }~~ A \in \Xi^{\lhd}(n,k). 
\end{eqnarray*}
In this case we can consider orbits in $\Xi(n,k)$ and in $\Xi^{\lhd}(n,k)$ under $\sigma$, 
while the parity may change.  \\
\indent
It was shown in [8] that to each permutation $A$ of $[n]$ there corresponds a smallest positive 
integer $\pi(A)$ such that $\sigma^{\pi(A)}A = A$, which is called the 
{\it period} of $A$. Its trace 
\[
           \{\sigma A, \sigma^2 A, \ldots, \sigma^{\pi(A)}A=A \} 
\]
is called the {\it orbit} of $A$. Also according to [8], if $A \in E^{\lhd}(n,k)$
the period $\pi(A)$ satisfies the relation 
\begin{eqnarray}
     \pi(A) =  (n-k)\gcd(n,\pi(A)).
\end{eqnarray}
\indent
From (i) and (ii) we see that orbits of permutations of $E^{\lhd}(n,k)$ 
contain those of the form $a_1a_2 \cdots a_{n-1}n$, 
where $a_1a_2 \cdots a_{n-1}$ are permutations of $[n-1]$. 
We call these permutations {\it canonical} in $E^{\lhd}(n,k)$ and in $\Xi^{\lhd}(n,k)$. 
Canonical permutations turn out useful in counting orbits in the following sections.\\
\begin{center}
{\bf 4. Proof of {\boldmath $D_{n,k}=R_{n,k}$}}
\end{center}
In this section we exhibit close relationships between PAPs and signed Eulerian numbers.
The first theorem states that the latter half of (2) holds for $R_{n,k}$. \\
\\
{\bf Theorem 4.1} {\it If $n$ is even and $1 \le k \le n-2$, the following relation holds:}
\begin{eqnarray*}
     R_{n,k} = R_{n-1,k-1} - R_{n-1,k}.
\end{eqnarray*}
\noindent
{\bf Proof.} 
First observe that it is possible to consider orbits in the set $\Xi(n,k)$ under $\sigma$, 
as was shown in the previous section, when $n$ is even.
For the proof we divide all permutations $a_1a_2 \cdots a_n$ of $\Xi_{\rm e}(n,k)$ 
into the three types by the position of $n$: 
\vspace{-0.5cm}
\begin{itemize}
\item[(i)]  $a_i = n$ for some $i$ $(2 \le i \le n-1)$; 
\item[(ii)]   $a_n = n$;
\item[(iii)] $a_1 = n$.
\end{itemize}
\vspace{-0.5cm}
On the other hand, we divide all permutations $a_1a_2 \cdots a_n$ of 
$\Xi_{\rm o}(n,k)$ into the following three types by the position of one: 
\vspace{-0.5cm}
\begin{itemize}
\item[(iv)]  $a_i = 1$ for some $i$ $(2 \le i \le n-1)$;
\item[(v)] $a_1 = 1$;
\item[(vi)] $a_n = 1$. 
\end{itemize}
\vspace{-0.5cm}
By evaluating the cardinality of the set of permutations for each type,
we will deduce $R_{n,k} = R_{n-1,k-1} - R_{n-1,k}$.\\
\indent
Let $A = a_1a_2 \cdots a_n$ be an even permutation of type (i). 
Then $\sigma A$ is odd one of type (iv), because by (4)
the difference of the numbers of inversions between $A$ and $\sigma A$ is
$|n+1-2i|$ and it is odd by assumption. Since $\sigma$ is a bijection, we see that with each
permutation $A$ of type (i) in $\Xi_{\rm e}(n,k)$ associates only one permutation 
of type (iv) in $\Xi_{\rm o}(n,k)$. Therefore, both types are irrelevant to 
the difference $R_{n,k} = |\Xi_{\rm e}(n,k)|-|\Xi_{\rm o}(n,k)| = P_{n,k} - Q_{n,k}$. \\
\indent
Let $A = a_1a_2 \cdots a_{n-1}n$ be of type (ii), where $a_1a_2 \cdots a_{n-1}$ is a
permutation of  $[n-1]$. Since $A$ is an even 
permutation, so is the permutation $a_1a_2 \cdots a_{n-1}$. 
Hence the total number of $\Xi_{\rm e}(n,k)$
of type (ii) is $P_{n-1,k-1}$. \\
\indent
Let $A = na_1a_2 \cdots a_{n-1}$ be of type (iii). Since $n-1$ is odd and $A$ is an even 
permutation, we see that $a_1a_2 \cdots a_{n-1}$ is an odd
permutation of  $[n-1]$. Therefore, the total number of $\Xi_{\rm e}(n,k)$
of type (iii) is $Q_{n-1,k}$. \\
\indent 
Similar arguments can be applied to types (v) and (vi). The total numbers 
of $\Xi_{\rm o}(n,k)$ of type (v) and type (vi) are equal to $Q_{n-1,k-1}$ and $P_{n-1,k}$, respectively. \\
\indent
Therefore, we conclude that the difference $R_{n,k}$ is described by
\[
    R_{n,k} = (P_{n-1,k-1} + Q_{n-1,k}) - (Q_{n-1,k-1} + P_{n-1,k}) =  R_{n-1,k-1} - R_{n-1,k},
\]
which is the required relation. 
~$\Box$ \\
\\
\indent
In order to show $D_{n,k}=R_{n,k}$, we first consider the case of odd $n$.\\
\\
{\bf Lemma 4.2.} {\it
If $n$ is odd and $0 \le k \le n-1$, then
\begin{eqnarray*}
       D_{n,k}=R_{n,k}.
\end{eqnarray*}
}
{\bf Proof.}
When $k=0$ or $n-1$, it is trivial, since such a permutation is a PAP. So assume $1 \le k \le n-2$.
With each permutation $A=a_1a_2\cdots a_n \in E(n,k)$ we associate canonical one
\[
B=a_1a_2\cdots a_n(n+1) \in E^{\lhd}(n+1,k+1).
\]
Obviously, we have ${\rm inv}(A) = {\rm inv}(B)$.
Using the operator $\sigma$, let us define another operator $\tau$ on the set of canonical permutations of 
$E^{\lhd}(n+1,k+1)$ by
\begin{eqnarray}
\tau B = \sigma^{n+1-a_n}B = b_1b_2\cdots b_n(n+1).
\end{eqnarray}
Notice that, applying $\tau$, the entry $a_n$ of $B$ changes into $n+1$ at the right end of $\tau B$ and
hence $\tau B$ is also canonical.
Here $b_1 =  n+1-a_n$ and $b_i = a_{i-1} + (n+1-a_n)$ $\pmod {n+1}$ for $i$ $(2 \le i \le n)$.
The operator $\tau$ also preserves the number of ascents of all permutations in $E^{\lhd}(n+1,k+1)$. \\
\indent
We examine the parity of $\tau B$. 
First, remark that ${\rm inv}(\sigma B) = {\rm inv}(B)$, for $B$ is canonical.
Each additional application of $\sigma$ 
changes the parity of permutations as long as $n+1$ lies in the interior of permutations, as shown
in Section 3, for $n+1$ is even.
Therefore, after $a_n$ of $B$ becomes $n+1$ by the application of $\sigma^{n+1-a_n}$ at the
right end of  $\tau B$, the parity of $B$ has changed $n-a_n$ times. Hence,
when $a_n$ is even, the parity of $\tau B$ is different from that of $B$, and 
when $a_n$ is odd, the parity of $\tau B$ is the same as that of $B$.\\
\indent
Suppose $A=a_1a_2\cdots a_n \in \Xi(n,k)$. Then $a_n$ is odd and $n-a_n$ is even, since $n$ is odd. 
Hence the corresponding $B$ to $A$ and $\tau B$ are PAPs, and the parity of $A$, $B$ and 
$\tau B$ is the same. We can continue applying this argument 
to $\tau B$,  $\tau^2 B$ and so on. Consequently, the orbit of $B$ under $\tau$
is a subset of canonical permutations of $E^{\lhd}(n+1,k+1)$ with the same parity. 
If a PAP is $A=14523 \in \Xi(5, 3)$, for example, then $B=145236 \in E^{\lhd}(6, 4)$, 
$\tau B=341256$ and $\tau^2B=145236$. The orbit of $B$ is $\{ B, \tau B\}$ under $\tau$. \\
\indent
Suppose, on the contrary, $A \in E(n,k)$ is not a PAP. Then neither is the corresponding $B$ and hence,
for some $\ell$, the parity of $\tau^{\ell}B$ turns opposite to that of $B$. 
The correspondence between $A$ and $B$ is a bijection, and $B$ changes the parity under $\tau$
among permutations that are not parity-alternate in $E^{\lhd}(n+1,k+1)$. More precisely,
if even (or odd) $A$ is not a PAP, the parity of its corresponding $B$ turns odd (or even) under $\tau$.
So we conclude that, among permutations that are not PAPs, 
the number of even permutations must be equal to that of odd ones. \\
\indent
Therefore, the difference $D_{n,k} = |E_{\rm e}(n,k)|- |E_{\rm o}(n,k)|$ depends 
only on the numbers of even PAPs and odd PAPs. This follows from the fact that
$B$'s for parity-alternate $A$'s do not change the parity under $\tau$.
Hence we obtain 
\[
D_{n,k} = P_{n,k} - Q_{n,k} = R_{n,k}, 
\]
which completes the proof.  ~$\Box$  \\
\\
\indent
The main result concerning signed Eulerian numbers in relation to PAPs
immediately follows from Theorem 4.1 and Lemma 4.2. \\
\\
{\bf Theorem 4.3.} {\it For all $n$ and $k$ {\rm (}$0 \le k \le n-1${\rm)} it follows that
\begin{eqnarray}
     R_{n,k} = D_{n,k}.
\end{eqnarray}
The recurrence relation for $R_{n,k}$ is given by
\begin{eqnarray*}
     R_{n,k} = 
        \left\{\begin{array}{ll}
                    (n-k)R_{n-1,k-1}+ (k+1)R_{n-1,k}, & \mbox{if  $n$ is odd}, \\
                    R_{n-1,k-1} - R_{n-1,k}, &  \mbox{if  $n$ is even}. \\
                    \end{array}  \right.
\end{eqnarray*}
}
\\
{\bf Proof.} For (7) it suffices to prove the case of even $n$. From (2) and Lemma 4.2 we have
\[
D_{n,k} = D_{n-1,k-1} - D_{n-1,k}=R_{n-1,k-1}-R_{n-1,k},
\]
which is equal to $R_{n,k}$ by Theorem 4.1. The recurrence relation for $R_{n,k}$ directly follows from (2).
 ~$\Box$  \\
\begin{center}
{\bf 5. Applications of {\boldmath $\tau$}}
\end{center}
In this section $n$ is assumed to be an even integer. Among canonical permutations in $\Xi^{\lhd}(n,k)$
we define the operator as in (6). That is, for $A = a_1\cdots a_{n-1}n \in \Xi^{\lhd}(n,k)$ let us put
\[
\tau A=\sigma^{n-a_{n-1}}(a_1\cdots a_{n-1}n) = b_1 \cdots b_{n-1}n,
\]
where $b_1 = n-a_{n-1}$ and $b_i=a_{i-1}+n-a_{n-1} \pmod{n}$ for $i$ ($2 \le i \le n-1$).
Since $n$ is even, $\tau$ preserves the number of
ascents and the parity of canonical PAPs, as was seen in the proof of Lemma 4.2. 
Therefore, we can consider orbits under $\tau$ among canonical permutations of 
$\Xi^{\lhd}_{\rm e}(n,k)$ and of $\Xi^{\lhd}_{\rm o}(n,k)$.
\\
\indent
It is easy to see that $\tau^n A = A$, since the last entry $n$ of canonical PAPs
moves to the left end by each application of $\tau$ and 
hence all entries return to the original positions after $\tau^n$.  
Therefore, the period of orbit of $A$ under $\tau$ is a divisor $d$ of $n$. 
This enables us to count orbits of period of each divisor $d$ of $n$. 
Let us denote by $\alpha^k_d$ the number of orbits of period $d$ among even canonical PAPs of 
$\Xi_{\rm e}^{\lhd}(n,k)$ and by $\beta^k_d$  that among odd 
canonical ones in $\Xi_{\rm o}^{\lhd}(n,k)$. A condition for these numbers to be zero
is deduced. 
Making use of this, we will deduce divisibility properties for $P_{n,k}$, 
$Q_{n,k}$ and $R_{n,k}$ by prime powers. Finally, when $n$ is even, a relation similar to (1) is
derived for $S_{n,k}$. \\
\\
{\bf Theorem 5.1.} {\it
Let $n$ be an even integer and let $k$ be an integer such that $1 \le k \le n-1$.
If a divisor $d$ of $n$ satisfies $\gcd(k, n/d) > 1,$ 
then $\alpha_d^k = 0$ and $\beta_d^k = 0$. 
}\\
\\
\noindent
{\bf Proof.} 
First we show that the period of the orbit of a canonical $A \in \Xi^{\lhd}(n,k)$ under $\tau$ is $d$ 
if and only if $d(n-k)$ is the period of the orbit of $A$ under $\sigma$. 
Suppose the period of $A$ under $\sigma$ is $d(n-k)$. 
By [8, Theorem 1 and Corollary 2] we see that the orbit of $A$ under $\sigma$ repeates itself $n/d$ times in
\[
 \{ \sigma A, \sigma^2 A, \ldots, \sigma^{n(n-k)}A=A \},
\]
where we find $n$ canonical permutations. Therefore, the orbit under $\sigma$ 
\begin{eqnarray}
\{ \sigma A, \sigma^2 A, \ldots, \sigma^{d(n-k)}A=A \}
\end{eqnarray}
has $d$ canonical ones. Hence the period under $\tau$ is equal to or less than $d$. \\
\indent
Next suppose the period of such a canonical $A$ under $\tau$ is $d$. 
This means that there are $d$ canonical permutations in (8). 
Then we have $\sigma^{d(n-k)}A=A$.
Hence the orbit under $\sigma$ is equal to or less than $d(n-k)$. It follows that 
the period of a canonical $A \in \Xi^{\lhd}(n,k)$ under $\tau$ is $d$
if and only if $d(n-k)$ is the period of $A$ under $\sigma$. \\
\indent
Assume that $A$ is an even canonical PAP of
$\Xi_{\rm e}^{\lhd}(n,k)$. If the period $\pi(A)$ under $\sigma$ is equal to $d(n-k)$ by (5), then it satisfies
$\pi(A) = (n-k)\gcd(n, \pi(A))$. Putting 
$d = \gcd(n, \pi(A))$, we have $\pi(A) = d(n-k)$ and $d = \gcd(n, d(n-k))$, 
which implies $\gcd(n-k, n/d) = 1$ or $\gcd(k, n/d) = 1$. 
Consequently, we see that there exist no permutations of 
period $d(n-k)$, i.e., $\alpha_d^k = 0$, if a divisor $d$ of $n$ satisfies 
$\gcd(k, n/d) > 1$. Similar arguments can be employed for proving $\beta^k_d = 0$, by
dealing with permutations in $\Xi_{\rm o}^{\lhd}(n,k)$.  ~$\Box$\\
\\
\noindent
{\bf Theorem 5.2.} {\it
Let $p$ be a prime and an even integer $n$ be divisible by $p^m$ for a positive integer $m$. 
If $k$ is divisible by $p$, then $S_{n -1, k-1}$, $P_{n -1, k-1}$, $Q_{n -1, k-1}$, $R_{n -1, k-1}$ 
and $D_{n -1, k-1}$
are also divisible by $p^m$.
}\\
\\
\noindent
{\bf Proof.}  
Without loss of generality we can assume 
that $m$ is the largest integer for which $p^m$ divides $n$. 
Suppose $k$ is a multiple of $p$. \\
\indent
First consider even PAPs in $\Xi_{\rm e}^{\lhd}(n,k)$. 
It suffices to deal only with canonical ones in counting orbits.
If $A=a_1a_2\cdots a_{n-1}n \in \Xi_{\rm e}^{\lhd}(n,k)$, we see that
$a_1a_2 \cdots a_{n-1} \in \Xi_{\rm e}(n-1,k-1)$. 
Therefore, there exist $P_{n-1,k-1}$ canonical permutations in $\Xi_{\rm e}^{\lhd}(n,k)$. 
Since there exist $\alpha_d^{k}$ orbits of period $d$ of $n$ under $\tau$,
classifying all canonical permutations of $\Xi_{\rm e}^{\lhd}(n,k)$ into orbits leads 
us to the following:
\begin{eqnarray}
P_{n-1,k-1} = \sum_{d|n} d \alpha_d^k. 
\end{eqnarray}
In Theorem 5.1 we have seen that $\alpha_d^k = 0$ for a divisor $d$ of $n$ 
such that $\gcd(k, n/d)>1$.
On the other hand, it follows that a divisor $d$ for which $\gcd(k, n/d) =1$ 
must be a multiple of $p^m$, since $k$ is a multiple of $p$. Therefore, 
$P_{n -1, k -1}$ is divisible by $p^m$. 
Similarly, dealing with permutations in $\Xi_{\rm o}^{\lhd}(n,k)$,
the same property holds for $Q_{n -1, k -1}$ and we get
\begin{eqnarray}
Q_{n-1,k-1} = \sum_{d|n} d \beta_d^k. 
\end{eqnarray}
From the relations
\begin{eqnarray*}
 & S_{n -1, k -1} & = ~P_{n -1, k -1} + Q_{n -1, k -1},\\
 & D_{n -1, k -1} & =~ R_{n -1, k -1} =~ P_{n -1, k -1} - Q_{n -1, k -1}, 
\end{eqnarray*}
it follows that $S_{n-1, k-1}$,
$R_{n-1, k-1}$ and $D_{n-1,k-1}$ are also divisible by $p^m$.  ~$\Box$ \\
\\
\indent
This theorem enables us to prove the following divisibility property. For classical Euleian numbers
it was proved in [9, Theorem 7]. \\
\\
\noindent
{\bf Corollary 5.3.} {\it
Let $p$ be an odd prime and an integer $n$ be divisible by $p^m$ for a positive integer $m$.
If $k$ is divisible by $p$, then $B_{n -1, k-1}$ and $C_{n -1, k-1}$ are also divisible by $p^m$.
}\\
\\
\noindent
{\bf Proof.}  When $n$ is odd, it was proved in [10, Corollary 5.2].
So assume $n$ is even. By [9, Theorem 7] Eulerian numbers $A_{n-1,k-1}$ is divisible by $p^m$. 
Since $p$ is odd and
\[
B_{n -1, k-1} = (A_{n -1, k-1} + D_{n -1, k-1})/2, ~~ C_{n -1, k-1} = (A_{n -1, k-1} - D_{n -1, k-1})/2,
\]
we see that these two numbers are divisible by $p^m$ from Theorem 5.2.
~$\Box$\\
\\
\noindent
{\bf Corollary 5.4.} {\it
When $n$ is even, $S_{n,k} =  (n-k)S_{n-1,k-1} + (k+1)S_{n-1,k}$  holds.
}
\\
\\
\noindent
{\bf Proof.}  
Remarking (8), all permutations in $\Xi_{\rm e}^{\lhd}(n,k)$ are classified by $\alpha_d^k$ orbits of periods  
$d(n-k)$. Then, using (9), we have
\[
|\Xi_{\rm e}^{\lhd}(n,k)|= \sum_{d|n} d(n-k) \alpha_d^k = (n-k) \sum_{d|n} d \alpha_d^k =
  (n-k)P_{n-1,k-1}.
\]
Similarly, using (10), we have $|\Xi_{\rm o}^{\lhd}(n,k)|= (n-k)Q_{n-1,k-1}$. Since 
$|\Xi^{\lhd}(n,k)|$ is the sum of both numbers, we get
\begin{eqnarray}
|\Xi^{\lhd}(n,k)|= (n-k)(P_{n-1,k-1} + Q_{n-1,k-1}) = (n-k)S_{n-1,k-1}.
\end{eqnarray}
On the other hand, all permutations $A = a_1a_2\cdots a_n \in \Xi(n,k)$ with $a_1 > a_n$
are tranformed into those in $\Xi^{\lhd}(n,n-k-1)$ by the reflection $A^{\ast} = a_n\cdots a_2a_1$.
Therefore, from (11) its cardinality turns equal to $(n-(n-k-1))S_{n-1,n-k-1-1} = (k+1)S_{n-1, n-k-2}$,
which is just $(k+1)S_{n-1,k}$ by the symmetry property (3).
Consequently, $S_{n,k}=|\Xi(n,k)|$ can be written as the sum of $(n-k)S_{n-1,k-1}$ and $(k+1)S_{n-1,k}$.
This completes the proof.  ~$\Box$\\

\begin{center}
\section*{\normalsize References}
\end{center}
\begin{itemize}
\item[{[1]}] J. D\'esarm\'enien and D. Foata, 
The signed Eulerian numbers. {\it Discrete Math.} 99: 49-58, 1992.
\item[{[2]}]  D. Foata and M.-P.  Sch\"utzenberger, 
{\it Th\'eorie G\'eom\'etrique des Polyn\^omes Eul\'eriens. } Lecture Notes in
Mathematics, Vol. 138, Springer-Verlag, Berlin, 1970.
\item[{[3]}]  R. L. Graham, D. E. Knuth and  O. Patashnik, {\it Concrete Mathematics.}
Addison-Wesley, Reading, 1989.
\item[{[4]}]  A. Kerber, {\it Algebraic Combinatorics Via Finite Group Actions.}
BI-Wissenschafts- verlag, Mannheim, 1991. 
\item[{[5]}]  D. E.  Knuth, {\it The Art of Computer Programming, Vol. 3, Sorting and
Searching.} Addison-Wesley, Reading, 1973.
\item[{[6]}]  L. Lesieur and J.-L. Nicolas,  On the Eulerian numbers
    $M_n =\max_{1\le k \le n}A(n,k)$.  {\it Europ. J. Combin. } 13: 379-399, 1992.
\item[{[7]}]  J.-L. Loday,  Op\'erations sur l'homologie cyclique des alg\`ebres
commutatives. {\it Invent. Math. } 96: 205-230, 1989. 
\item[{[8]}]  S. Tanimoto,  An operator on permutations and its 
application to Eulerian numbers. {\it Europ. J. Combin. } 22: 569-576, 2001.
\item[{[9]}]  S. Tanimoto,  A study of Eulerian numbers by means of an operator on permutations.
{\it Europ. J. Combin. } 24: 33-43, 2003.
\item[{[10]}] S. Tanimoto,  A study of Eulerian numbers for permutations in the alternating group.
{\it Electronic J. of Comb. Number Theory } 6:  no.31, 2006. arXiv:math. CO/0602263.
\end{itemize}

\end{document}